\documentclass[12pt]{article}
\usepackage{amsmath,amsfonts,latexsym,amsthm,amssymb}
\numberwithin{equation}{section}
\DeclareMathOperator{\card}{card}
\DeclareMathOperator{\Tr}{Tr}
\DeclareMathOperator{\supp}{supp}
\DeclareMathOperator{\ord}{ord}
\DeclareMathOperator{\meas}{meas}
\DeclareMathOperator{\Dim}{Dim}
\DeclareMathOperator{\diam}{diam}

\newcommand{\K}{\overline{K}}
\begin{document}
\newtheorem{prop}{Proposition}
\newtheorem{lem}{Lemma}
\newtheorem{teo}{Theorem}
\newtheorem{cor}{Corollary}
\pagestyle{plain}
\begin{titlepage}

\vspace*{3cm} \noindent
{\bf {\Large Hausdorff Measure for a Stable-Like
Process \\ over an Infinite Extension of a Local
Field}}\footnote[1]{This research was supported in part by CRDF
Grant UM1-2090.}

\vspace*{1cm} \noindent
{\bf Anatoly N. Kochubei}\footnote[2]{Institute of Mathematics,
National Academy of Sciences of Ukraine, Tereshchenkivska 3, Kiev,
01601 Ukraine. E-mail: ank@ank.kiev.ua}

\vspace*{2cm}
Running head:\quad  ``Stable-Like
Process over an Infinite Extension of a Local Field''
\end{titlepage}
\vspace*{8cm}
\begin{abstract}
We consider an infinite extension $K$ of a local field of zero
characteristic which is a union of an increasing sequence of finite
extensions. $K$ is equipped with an inductive limit topology;
its conjugate $\K$ is a completion of $K$ with respect to
a topology given by certain explicitly written seminorms. The
semigroup of measures, which defines a stable-like process $X(t)$ on
$\K$, is concentrated on a compact subgroup $S\subset \K$. We
study properties of the process $X_S(t)$, a part of $X(t)$ in
$S$. It is shown that the Hausdorff and
packing dimensions of the image of an interval equal 0 almost surely.
In the case of tamely ramified extensions a correct Hausdorff
measure for this set is found.
\end{abstract}
\vspace{2cm}
\newpage
\section{INTRODUCTION}

Let $k$ be a non-Archimedean local field (= non-discrete totally
disconnected locally compact topological field) with
characteristic zero. Consider a strictly increasing sequence of
its finite algebraic extensions
\begin{equation}
k=K_1\subset K_2\subset \ldots \subset K_n\subset \ldots .
\end{equation}
The infinite extension
$$
K=\bigcup\limits_{n=1}^\infty K_n
$$
may be considered as a topological vector space over $k$ with the
inductive limit topology. Let $\K$ be its strong dual. Evidently,
$\K$ is not locally compact.

Within the non-Archimedean version of infinite-dimensional
analysis developed by the \linebreak author$^{(8-10)}$, an analog of the
symmetric $\alpha$-stable process was constructed on $\K$. Just
as for stable processes on local fields (see e.g. Refs. 1, 5, 7, 10,
14, 17), this process $X(t)$ is defined for any $\alpha >0$. Some of its
properties are similar to those of the classical stable processes
or processes on local fields while others are different. In
particular, $X(t)$ has an invariant measure $\mu$ which is
Gaussian in the sense of Evans$^{(2)}$; the transition
probabilities of $X(t)$ are not absolutely continuous with
respect to $\mu$.

Note that both $\mu$ and the convolution semigroup of measures
$\pi (t,dx)$, which defines $X(t)$, are concentrated on a compact
subgroup $S\subset \K$, and $\mu$ coincides with the normalized
Haar measure on $S$ (for the details see Section 2 below). Thus
an essential information on the process $X(t)$ is contained in
the properties of its part $X_S(t)$ in $S$.

In order to study sample path properties of $X_S(t)$, we can use
the results by Evans$^{(3)}$ who investigated L\'evy processes on
a general Vilenkin group (a non-discrete locally compact totally
disconnected Abelian topological group). The topology in a
Vilenkin group is determined by a descending chain of compact
open subgroups. This chain is not unique, and as soon as we
manage to write such a chain $\{S_n\}$ explicitly for our case
(Section 3) and compute the L\'evy measure of $S\setminus S_n$
(Section 4), the general theorems by Evans$^{(3)}$ yield
immediately the asymptotics of the first exit time $\pi (n)$ of
$X_S(t)$ out of the subgroup $S_n$, and an information on the
local behavior of sample paths. We also prove that both the
Hausdorff and packing dimensions of a sample path of $X_S(t)$
equal 0 almost surely, which is quite different both from the
classical case and the case of a local field considered recently
by Albeverio and Zhao$^{(1)}$.

The last result shows the importance of finding, for our
situation, a correct Hausdorff measure. However this problem is
more complicated. In order to use the appropriate theorem by
Evans$^{(3)}$, we have to know that
\begin{equation}
\liminf\limits_{n\to \infty }Q(n,N)>0
\end{equation}
where
$$
Q(n,N)=\mathbf P\left\{ X_S(t)\notin S_n\ \forall \ t\in [\pi
(n),\pi (N))\right\} ,\quad n>N.
$$
Evans proved (1.2) for processes with locally spherically
symmetric L\'evy measures. This condition is not fulfilled in our
case, and we give (Section 5) a direct proof of (1.2) under the
assumption that all the extensions in (1.1) are tamely ramified.
Such an assumption is often made in algebraic number theory
because the algebraic structure of tamely ramified extensions is
more or less transparent while general extensions may behave
quite wildly.

The author is grateful to the referee for very helpful comments
and suggestions.

\section{PRELIMINARIES}

{\bf 2.1.} Let us recall some of the main constructions and
results from Refs. 8, 9. For the basics on local fields see Refs.
4, 12, 13, 16. We also use some standard material regarding field
extensions which can be found in algebra textbooks (Refs. 11, 15).

Consider the sequence of extensions (1.1). For each
$n=1,2,\ldots$, we define a mapping $T_n:\ K\to K_n$ as follows.
If $x\in K_\nu$, $\nu >n$, put
$$
T_n(x)=\frac{m_n}{m_\nu}\Tr_{K_\nu /K_n}(x)
$$
where $m_n$ is the degree of the extension $K_n/k$, $\Tr_{K_\nu
/K_n}:\ K_\nu \to K_n$ is the trace mapping. If $x\in K_n$, then,
by definition, $T_n(x)=x$. The mapping $T_n$ is well-defined, and
$T_n\circ T_\nu =T_n$ for $\nu >n$. We shall write $T$ instead of
$T_1$.

The strong dual space $\K$ can be identified with the projective
limit of the sequence $\{K_n\}$ with respect to the mappings
$\{T_n\}$, that is with the subset of the direct product
$\prod\limits_{n=1}^\infty K_n$ consisting of those $x=(x_1,\ldots
,x_n,\ldots )$, $x_n\in K_n$, for which $x_n=T_n(x_\nu )$ if $\nu
>n$. The topology in $\K$ is defined by seminorms
$$
\|x\|_n=\|x_n\|,\quad n=1,2,\ldots ,
$$
where $\|\cdot \|$ is the extension onto $K$ of the normalized
absolute value $|\cdot |_1$ defined on $k$. If $x\in K_n$ then
$\|x\|=|x|_n^{1/m_n}$ where $|x|_n$ is the normalized absolute
value on $K_n$.

The pairing between $K$ and $\K$ is defined as
$$
\langle x,y\rangle =T(xy_n)
$$
where $x\in K_n\subset K$, $y=(y_1,\ldots ,y_n,\ldots )\in \K$,
$y_n\in K_n$. Identifying an element $x\in K$ with $(x_1,\ldots
,x_n,\ldots )\in \K$ where $x_n=T_n(x)$, we can view $K$ as a
dense subset of $\K$. The mappings $T_n$ can be extended to
linear continuous mappings from $\K$ to $K_n$, by setting
$T_n(x)=x_n$ for any $x=(x_1,\ldots ,x_n,\ldots )\in \K$.

Consider a function on $K$ of the form
$$
\Omega (x)=\left\{ \begin{array}{rl}
1 ,& \mbox{if } \|x\|\le 1\\
0, & \mbox{if } \|x\|>1
\end{array}\right.
$$
$\Omega$ is continuous and positive definite on $K$. That results
in the existence of a probability measure $\mu$ on the Borel
$\sigma$-algebra $\mathcal B(\K)$, such that
$$
\Omega (a)=\int\limits_{\K}\chi (\langle a,x\rangle )\,\mu
(dx),\ \ a\in K,
$$
where $\chi$ is a rank zero additive character on $k$. The
measure $\mu$ is Gaussian in the sense of Ref. 2. It is
concentrated on the compact additive subgroup
$$
S=\left\{ x\in \K:\ \|T_n(x)\|\le q_n^{d_n/m_n}\|m_n\|,\
n=1,2,\ldots \right\}
$$
where $q_n$ is the residue field cardinality for the field $K_n$,
$d_n$ is the exponent of the different of the extension $K_n/k$.

We shall call $S$ {\it the support subgroup} of $\K$. The
restriction of $\mu$ to $S$ coincides with the normalized Haar
measure on $S$. If $f$ is a ``cylindrical'' function on $\K$ of
the form $f(x)=\varphi (T_n(x))$, $x\in \K$, where $\varphi$ is a
locally integrable complex-valued function on $K_n$, then
\begin{equation}
\int\limits_{\K}f(x)\mu (dx)=q_n^{-d_n}\|m_n\|^{-
m_n}\int\limits_{z\in K_n:\ \|z\|\le q_n^{d_n/m_n}\|m_n\|}\varphi
(z)\,dz
\end{equation}
($dz$ is the Haar measure on $K_n$ normalized by the condition
$\int\limits_{|z|_n\le 1}dz=1$).

Denote
$$
\rho _\alpha (s,t)=\left\{ \begin{array}{rl}
e^{-ts^\alpha }, & \mbox{if } s>1 \\
1, & \mbox{if } s\le 1,
\end{array} \right.
$$
where $\alpha >0$, $t>0$. For each $t>0$ the function $\rho
_\alpha (\|\xi \|,t)$ is a continuous positive definite function
on $K$, and there exists a family $\pi (t,dx)$ of Radon
probability measures on $\mathcal B(K)$ (concentrated on $S$)
such that
\begin{equation}
\rho _\alpha (\|\xi \|,t)=\int\limits_{\K}\chi (\langle \xi ,x
\rangle )\pi (t,dx),\quad \xi \in K.
\end{equation}
This family is actually a semigroup of measures defining a L\'evy
process $X(t)$ on $\K$. Its part in $S$ will be denoted by
$X_S(t)$. We have an integral formula similar to (2.1): if
$f(x)=\varphi (T_n(x))$ then
\begin{equation}
\int _{\overline{K}}f(x)\pi (t,dx)
=\|m_n\|^{-m_n}\int \limits_{z\in K_n:\ \|z\|
\le q_n^{d_n/m_n}\|m_n\|}\Gamma _\alpha ^{(n)}(m_n^{-1}z,t)\varphi (z)\,dz
\end{equation}
where
$$
\Gamma _\alpha ^{(n)}(x,t)=q_n^{-d_n}\int\limits_{K_n}(\chi
\circ \mbox{Tr}_{K_n/k})(-x\xi )\rho _\alpha (\|\xi \|,t)\,d\xi .
$$

The generator of the process $X(t)$ is a hyper-singular integral
operator which resembles the fractional differentiation operator
$D^\alpha$ generating the symmetric stable process on a local
field. The generator can be expressed in terms of the L\'evy
measure $\Pi (dx)$ on $\mathcal B(\K\setminus \{0\})$ which
appears in the following formula of the L\'evy-Khinchin type: for
any $\lambda \in K$, $t>0$
$$
\mathbf E\chi (\langle \lambda ,X(t)\rangle )=\exp \left\{
t\int\limits_{\K}[\chi (\langle \lambda ,x\rangle )-1]\Pi
(dx)\right\} .
$$
We have the identity
\begin{equation}
\int\limits_{\K}[\chi (\langle \lambda ,x\rangle )-1]\Pi (dx)=
\left\{ \begin{array}{rl}
-\|\lambda \|^\alpha ,& \mbox{if } \|\lambda \|>1\\
0, & \mbox{if } \|\lambda \|\le 1
\end{array}\right.
\end{equation}

As it could be expected, the measure $\Pi$ is concentrated on
$S\setminus \{0\}$. Indeed, let
$$
S^{(n)}=\left\{ x\in \K:\ \|T_n(x)\|\le
q_n^{d_n/m_n}\|m_n\|\right\} ,\quad n=1,2,\ldots .
$$
Then (Ref. 8) $S^{(\nu )}\subset S^{(n)}$ for $\nu >n$, and
$S=\bigcap_{n=1}^\infty S^{(n)}$, so that $\K\setminus
S=\bigcup_{n=1}^\infty \left( \K\setminus S^{(n)}\right)$. It was
shown in Ref. 9 that for a function $f$ as in (2.1) or (2.3) with
$0\notin \supp \varphi$ we have
\begin{multline}
\int \limits _{\K}f(x)\Pi (dx) =
-q_n^{d_n\alpha /m_n}\frac{1-q_n^{\alpha /m_n}}{1-q_n^{-1-
\alpha /m_n}}\int \limits_{x\in K_n:\  |x|_n
\le q_n^{d_n}}\biggl[ |x|_n^{-1-\alpha /m_n}\biggr. \\ +\biggl.
\frac{1-q_n^{-1}}{q_n^{\alpha /m_n}-1}q_n^{-d_n(1+\alpha /m_n)}
\biggr] \varphi (m_nx)\,dx.
\end{multline}
It follows from (2.5) that $\Pi \left( \K\setminus S^{(n)}\right)
=0$ whence $\Pi \left( \K\setminus S\right) =0$.

We shall also use another analytic expression for the integral in
the left-hand side of (2.5):
\begin{equation}
\int \limits _{\K}f(x)\Pi (dx) =-\int\limits_{\eta \in K_n:\
\|\eta \|>1}\|\eta \|^\alpha w_n(\eta )\,d\eta
\end{equation}
where $w_n$ is the inverse Fourier transform of the function
$y\mapsto q_n^{-d_n/2}\varphi (m_ny)$, that is
$$
w_n(\eta )=q_n^{-d_n}\int\limits_{K_n}\chi \circ \Tr_{K_n
/k}(-\eta y)\varphi (m_ny)\,dy.
$$

Indeed, after an obvious change of variables we can write the
Fourier inversion as
$$
\varphi (z)=\int\limits_{K_n}(\chi \circ T)(\eta z)w_n(\eta
)\,d\eta ,\quad z\in K_n.
$$
Since by our assumption $\varphi (0)=0$, we have
$$
\int\limits_{K_n}w_n(\eta )\,d\eta =0,
$$
so that for any $x\in \K$
$$
\varphi (T_n(x))=\int\limits_{K_n}[(\chi \circ T)(\eta
T_n(x))-1]w_n(\eta )\,d\eta
=\int\limits_{K_n}[(\chi (\langle x,\eta \rangle )-1]w_n(\eta
)\,d\eta .
$$
Integrating with respect to $\Pi$ and using (2.4) we come to
(2.6).

\bigskip
{\bf 2.2.} Here we collect, for a reader's convenience, some
notions regarding ramification in extensions of local fields.

Let $K$ be a finite extension of a local field $L$ (this is
traditionally denoted $K/L$); we consider only fields of zero
characteristic. Denote by $O_K$, $O_L$ the corresponding rings of
integers, and by $P_K$, $P_L$ the prime ideals. We have
$\Tr_{K/L}:\ O_K\to O_L$. Moreover, there exists such an integer
$d$ ({\it the exponent of the different} for the extension $K/L$)
that $\Tr_{K/L}(x)\in O_L$ for $|x|_K\le q_K^d$ but
$\Tr_{K/L}(x_0)\notin O_L$ for some $x_0$ with
$|x_0|_K=q_K^{d+1}$. Here we furnish with appropriate subscripts
the objects related to $K$ or $L$ (the normalized absolute values
$|\cdot |_K$, $|\cdot |_L$, the residue field cardinalities
$q_K$, $q_L$ etc.). If $\pi_K$, $\pi_L$ are the
prime elements of $K$ and $L$, we have
$$
|\pi_K|_K=q_K^{-1},\quad |\pi_L|_L=q_L^{-1},\quad |\pi_L|_K=q_K^{-
e},
$$
where $e\ge 1$ is {\it the ramification index} of the extension $K/L$.

It is known that $d\ge e-1$, and $d=0$ if and only if $e=1$ (in
this case the extension $K/L$ is called {\it unramified}).
The extension $K/L$ is called {\it tamely ramified} if $d=e-1$.
$K/L$ is tamely ramified if and only if the characteristic of finite fields
$O_K/P_K$ and $O_L/P_L$ does not divide $e$.

Another important number attached to the extension $K/L$ is its
{\it index of inertia} $f\ge 1$ defined by the relation
$q_K=q_L^f$. The product $ef$ coincides with the degree $[K:L]$
of the extension. If $f=1$, the extension is called {\it totally
ramified}.

As an example, consider the case where $L$ is the field $\mathbb
Q_p$ of $p$-adic numbers, and $K$ is a quadratic extension, that
is an extension of the degree 2 obtained by adjoining an element
$\sqrt{\tau }$, where $\tau \in \mathbb Q_p$ is not a square of
an element of $\mathbb Q_p$, with natural algebraic operations
and the normalized absolute value
$$
\left| x_1+\sqrt{\tau }x_2\right|=\left| x_1^2-\tau
x_2^2\right|_{\mathbb Q_p},\quad x_1,x_2\in \mathbb Q_p.
$$
If $p\ne 2$, then there are 3 different extensions of the above
form, in which $|\tau |_{\mathbb Q_p}=1$ (the unramified
extension), or $\tau =p$, or $\tau =\varepsilon p$, $|\varepsilon
|_{\mathbb Q_p}=1$. The two latter extensions are totally and
tamely ramified ($e=2$ is prime to $p$). If $p=2$, the situation
is more complicated (see Proposition 5.12 in Ref. 12), there are
7 different quadratic extensions, one of which is unramified
while others are totally but not tamely ramified.

Let us return to the general situation.
Let $\widehat{O}_K\subset O_K$ be a complete system of
representatives of residue classes from $O_K/P_K$. Then any
nonzero element $x\in K$ is uniquely representable in the form of
the convergent series
\begin{equation}
x=\pi_K^{-n}(x_0+x_1\pi_K+x_2\pi_K^2+\cdots )
\end{equation}
where $n\in \mathbb Z$, $|x|_K=q_K^n$, $x_j\in \widehat{O}_K$,
$x_0\notin P_K$. If the extension $K/L$ is unramified then any
prime element of $L$ is simultaneously a prime element of $K$,
and in this case we may put $\pi_L$ instead of $\pi_K$ in (2.7).
On the other hand, if $K/L$ is totally ramified, then we may
identify $O_K/P_K$ with $O_L/P_L$, and use in (2.7) elements
$x_j\in \widehat{O}_L$. Note also that the Galois group of an
unramified extension is cyclic. Its generator is called {\it the
Frobenius automorphism} of the extension $K/L$.

If $K/L$ is an arbitrary finite extension, there exists an
intermediate extension $L'$ lying between $L$ and $K$, such that
$L'/L$ is unramified while $K/L'$ is totally ramified. In this
case the ramification index $e$ of $K/L$ equals the degree of
$K/L'$, the inertia index $f$ of $K/L$ coincides with the degree
of $L'/L$. The field $L'$ is called {\it the inertia subfield} of $K/L$.
If, in addition, $K/L$ is tamely ramified, then the prime elements
$\pi_K$ and $\pi_{L'}$ in $K$ and $L'$ can be chosen in such a way that
$\pi_K^e=\pi_{L'}$.

Dealing with the extensions (1.1) we shall denote the ramification index,
inertia index, and exponent of the different for the extension $K_\nu /K_n$,
$\nu >n$, by $e_{\nu n}$, $f_{\nu n}$, $d_{\nu n}$ respectively.
We shall write $e_\nu$, $f_\nu$, $d_\nu$ instead of $e_{\nu 1}$, $f_{\nu 1}$,
$d_{\nu 1}$. We shall denote by $O_n$ and $P_n$ respectively the
ring of integers and the prime ideal of $K_n$; $\widehat{O}_n$
will denote a complete system of representatives of residue classes from
$O_n/P_n$. Since $k$, as any local field of zero characteristic,
is a finite extension of the field $\mathbb Q_p$ of $p$-adic
numbers (here $p$ is the characteristic of the residue field), we
may assume without restricting generality that $k=\mathbb Q_p$.

\section{SUPPORT SUBGROUP AND ITS DUAL}

{\bf 3.1.} Our investigation of the group $S$ will be based on
the following auxiliary result.

Denote
$$
\Sigma_{n,N}=\left\{ y\in K_n:\ \|y\|\le q_n^{d_n/m_n-
N/f_n}\|m_n\|\right\},\quad n\ge 1,N\ge 0.
$$

\begin{lem}
If $\nu >n$, then $T_n$ maps $\Sigma_{\nu ,N}$ onto
$\Sigma_{n,N}$.
\end{lem}

{\it Proof}. Let $y\in \Sigma_{\nu ,N}$. Then
$$
|m_\nu^{-1}y|_\nu \le q_\nu^{d_\nu -Ne_\nu}.
$$
We have
$$
T_n(y)=\frac{m_n}{m_\nu }\Tr_{K_\nu
/K_n}(y)=m_n\Tr_{K_\nu /K_n}(m_\nu^{-1}y).
$$
Well-known properties of traces in local field extensions (Ref.
16, Chapter 8, Proposition 4) imply the inequality
$$
|T_n(y)|_n\le |m_n|_nq_n^l
$$
where $l\in \mathbb Z$ is determined from the inequality
$$
e_{\nu n}(l-1)<d_\nu -Ne_\nu -d_{\nu n}\le e_{\nu n}l.
$$
It is also known (Ref. 16, Chapter 8) that $d_\nu =e_{\nu n}d_n+d_{\nu n}$,
so that
\begin{equation}
e_{\nu n}(l-1)<e_{\nu n}d_n-Ne_\nu\le e_{\nu n}l.
\end{equation}

On the other hand (Ref. 4, Chapter II, (2.1)), $e_\nu =e_{\nu
n}e_n$, and we see from (3.1) that
$$
l-1<d_n-Ne_n\le l,
$$
whence $l=d_n-Ne_n$. It means that $T_n(y)\in \Sigma_{n,N}$ as
desired.

Conversely, if $z\in \Sigma_{n,N}$, that is $|m_n^{-1}z|_n\le
q_n^{d_n-Ne_n}$, it follows from the surjectivity property of the
trace (Ref. 16, Chapter 8, Proposition 4) that there exists such
an element $y'\in K_\nu$ that
$$
|y'|_\nu \le q_\nu^{d_\nu -Ne_\nu },\quad \Tr_{K_\nu
/K_n}(y')=m_n^{-1}z.
$$
Setting $y=m_\nu y'$ we find that
$$
z=\frac{m_n}{m_\nu }\Tr_{K_\nu /K_n}(y)=T_n(y),\quad |y|_\nu \le
|m_\nu |_\nu q_\nu^{d_\nu -Ne_\nu },
$$
which completes the proof. $\quad \blacksquare$

\medskip
Let us consider the sequence of subgroups
\begin{equation}
S=S_0\supset S_1\supset S_2\supset \ldots \supset \{0\},
\end{equation}
$$
S_n=\left\{ x\in S:\ \|T_n(x)\|\le q_n^{d_n/m_n-
n/f_n}\|m_n\|\right\},\quad n\ge 1.
$$
If $x\in S_n$ and $j<n$, then $T_j(x)=T_j(T_n(x))$, and by Lemma
1
$$
\|T_j(x)\|\le q_j^{d_j/m_j-n/f_j}\|m_j\|\le q_j^{d_j/m_j-
j/f_j}\|m_j\|,
$$
so that $x\in S_j$. We have shown that the subgroups $S_n$ indeed
form a descending chain. The same argument shows also that
$\bigcap\limits_{n=0}^\infty S_n=\{0\}$. The subgroups $S_n$ are
open and closed. It follows from Lemma 1 that the system of
subgroups $\{S_n\}$ forms a base of neighbourhoods of the origin
in $S$.

The quotient group $S/S_n$ is finite. Denote $M(n)=\ord S/S_n$.
It follows from the invariance of the measure $\mu$ that $\mu
(S_n)=[M(n)]^{-1}$. On the other hand, it follows from (2.1) that
$\mu (S_n)=q_1^{-nm_n}$. Thus
\begin{equation}
M(n)=q_1^{nm_n}.
\end{equation}

If we define $|x|$ for $x\in S$ by setting
$$
|x|=
\begin{cases}
1, & \text{if $x\notin S_1$}, \\
\left[ M(n)\right]^{-1}, & \text{if $x\in S_n\setminus S_{n+1}$,
$n=1,2,\ldots$},
\end{cases}
$$
and $|0|=0$, then $\Delta (x,y)=|x-y|$ is an ultrametric on $S$.

The descending chain (3.2) can be ``lengthened'' by including
intermediate subgroups so that the resulting chain would be such
that the quotient group of two consequtive subgroups is of a
prime order. This property was assumed in Ref. 3. However it will
be more convenient for us to use the chain (3.2). All the results
of Ref. 3 remain valid here.

\medskip
{\bf 3.2.} Let $S^*$ be the dual group of $S$. It was shown in
Ref. 9 that $S^*$ is isomorphic to $K/O$ where $O=\{\xi \in K:\
\|\xi \|\le 1\}=\bigcup\limits_{n=1}^\infty O_n$. Let
$\Xi_n\subset S^*$ be the annihilator of the subgroup $S_n$, that
is
$$
\Xi_n=\{ \xi +O:\ \xi \in K,\chi (\langle \xi ,x\rangle )=1\
\forall \ x\in S_n\}.
$$

\begin{prop}
If the extensions (1.1) are tamely ramified, then
$$
\Xi_n=\left\{ \xi +O:\ \xi \in K_n,\ |\xi |_n\le q_n^{ne_n}\right\}.
$$
\end{prop}

{\it Proof}. Assume that $k=\mathbb Q_p$. If $\xi \in K_n$,
$\ |\xi |_n\le q_n^{ne_n}$, then for any $x\in S_n$
$$
\left| m_n^{-1}\xi T_n(x)\right|_n=|\xi |_n\|m_n^{-
1}T_n(x)\|^{m_n}\le q_n^{ne_n}\cdot q_n^{d_n-ne_n}=q_n^{d_n},
$$
so that
\begin{equation}
|T(\xi T_n(x))|_1\le 1\quad \mbox{for all }x\in S_n,
\end{equation}
whence $\xi +O\in \Xi_n$.

Conversely, let $\xi +O\in \Xi_n$. Suppose first that $\xi \in
K_n$. Then (3.4) holds. Any element $z_n\in K_n$ with $\|z_n\|\le
q_n^{d_n/m_n-n/f_n}\|m_n\|$ can be ``lifted'' to an element $x\in
S_n$ such that $T_n(x)=z_n$. Indeed by Lemma 1 there exists
$z_{n+1}\in \Sigma_{n+1,n}\subset \Sigma_{n+1,0}$ such that
$T_n(z_{n+1})=z_n$. Then we find $z_{n+2}\in \Sigma_{n+2,0}$,
$T_{n+1}(z_{n+2})=z_{n+1}$. Repeating this and setting
$z_j=T_j(z_n)$ for $j<n$ we obtain $x=(z_1,\ldots ,z_n,\ldots
)\in S$ with $T_n(x)=z_n$. It is clear that $x\in S_n$.

If $|\xi |_n>q_n^{ne_n}$, and $\|T_n(x)\|=q_n^{d_n/m_n-
n/f_n}\|m_n\|$, we have $|\xi T_n(x)|_n>q_n^{d_n}$, and we can
choose $x\in S_n$ in such a way that (3.4) is violated. This
proves that $|\xi |_n\le q_n^{ne_n}$.

Now it remains to prove that any element of $\Xi_n$ can be
represented as $\xi +O$ with $\xi \in K_n$. Suppose that $\xi +O\in
\Xi_n$, $\xi \in K_l$, $l>n$. Denote by $K_n'$ the inertia
subfield of the extension $K_l/K_n$. Then $[K_n':K_n]=f_{ln}$,
$[K_l:K_n']=e_{ln}$, so that $[K_n':K_1]=m_nf_{ln}$, the
ramification index $e(K_n',K_1)=e_n$, and the inertia index
$f(K_n',K_1)=f_nf_{ln}$. Considering the extensions $K_n'\supset
K_n\supset K_1$ (see Section VIII-1 in Ref. 16) we find also the exponent
of the different $d(K_n',K_1)=e(K_n',K_n)d_n+d(K_n',K_n)=d_n$. The
residue field cardinality for $K_n'$ equals $q_n^{f_{ln}}$.

Define for $z\in K_l$
$$
T_n'(z)=\frac{[K_n':K_1]}{[K_l:K_1]}\Tr_{K_l/K_n'}(z)=
\frac{m_nf_{ln}}{m_l}\Tr_{K_l/K_n'}(z).
$$
Just as in Lemma 1, if $z\in \Sigma_{l,0}$, then
$$
\left\| T_n'(z)\right\| \le
q_n^{f_{ln}d(K_n',K_1)/[K_n':K_1]}\left\| [K_n':K_1]\right\|
=q_n^{d_n/m_n}\|m_nf_{ln}\|.
$$
For any $u\in K_n'$ $|u|_{K_n'}=\|u\|^{m_nf_{ln}}$ whence
\begin{equation}
\left| T_n'(z)\right|_{K_n'}\le
q_n^{d_nf_{ln}}|m_nf_{ln}|_{K_n'}.
\end{equation}

Let us consider $u=T_n'(z)$ as an element of the field $K_l$. We
have by (3.5)
$$
|u|_l=\|u\|^{m_l}=\left( \|u\|^{m_nf_{ln}}\right)
^{\frac{m_l}{m_nf_{ln}}}=|u|_{K_n'}^{\frac{m_l}{m_nf_{ln}}}\le
q_n^{\frac{d_nm_l}{m_n}}|m_nf_{ln}|_{K_n'}^{\frac{m_l}{m_nf_{ln}}}.
$$
We know that $q_l=q_n^{f_{ln}}$, $\frac{m_l}{m_n}=f_{ln}e_{ln}$,
$d_l=d_ne_{ln}+d_{ln}$, so that
$$
q_n^{\frac{d_nm_l}{m_n}}=q_l^{d_l-d_{ln}}
$$
and
$$
|m_nf_{ln}|_{K_n'}=|m_nf_{ln}|_l^{\frac{m_nf_{ln}}{m_l}}.
$$
Therefore
$$
|u|_l\le q_l^{d_l-d_{ln}}|m_l|_l\cdot |e_{ln}|_l^{-1}.
$$

Since the extension $K_l/K_n$ is tamely ramified, we find that
$|e_{ln}|_l=1$, $|u|_l\le q_l^{d_l}|m_l|_l$, so that $u\in
\Sigma_{l,0}$. This means in particular that for any $z\in
\Sigma_{l,0}$ the element $v=z-T_n'(z)$ also belongs to
$\Sigma_{l,0}$. Lifting $v$ to an element $y\in S$ we see that
$$
T_n(y)=T_n\circ T_n'(y)=T_n\circ T_n'(v)=0.
$$
Thus $y\in S_n$ and $|T(\xi (z-T_n'(z)))|_1\le 1$.

Now we have
\begin{multline*}
T(\xi z)-T(\xi T_n'(z))=T(\xi z)-T\circ T_n'(\xi T_n'(z))=
T(\xi z)-T(T_n'(z)T_n'(\xi ))\\
=T(\xi z)-T\circ T_n'(zT_n'(\xi ))=T(z(\xi -T_n'(\xi)))
\end{multline*}
and
$$
|T(z(\xi -T_n'(\xi)))|_1\le 1
$$
for any $z\in \Sigma_{l,0}$, which implies the inequality
\begin{equation}
|\xi -T_n'(\xi)|_l\le 1
\end{equation}
(since the annihilator of $\Sigma_{l,0}$ in $K_l=K_l^*$ is $O_l$; see
Ref. 9).

It follows from (3.6) that $\xi \in K_n'+O$. Let $g$ be the Frobenius
automorphism of the unramified extension $K_n'/K_n$. For an arbitrary
$z\in K_n'$ with $\|z\|\le q_n^{d_n/m_n}\|m_nf_{ln}\|$ (the set similar
to $\Sigma_{l,0}$ for the field $K_n'$) we consider the element
$z-g^{-1}z$ and its lifting $x_g\in S$. Then
$$
T_n(x_g)=T_n(T_n'(x_g))=T_n(z-g^{-1}z)=0
$$
due to the invariance of traces with respect to elements of the Galois
group. Thus $x_g\in S_n$ and
$$
|T(\xi T_n'(x_g))|_1=|T(\xi (z-g^{-1}z))|_1\le 1.
$$

Since
$$
T(\xi (z-g^{-1}z))=T(\xi z)-T(\xi \cdot (g^{-1}z))
=T(\xi z)-T(g^{-1}(g\xi \cdot z))=T(z(\xi -g\xi )),
$$
we see as above that $\xi -g\xi \in O$. Let us write the canonical
representation
$$
\xi =\pi_n^{-r}\left( \sigma_0+\sigma_1\pi_n+\cdots
+\sigma_{r-1}\pi_n^{r-1}+\cdots \right)
$$
where $\pi_n$ is a prime element both for $K_n$ and $K_n'$, $\sigma_j$ are
representatives of residue classes from $O_n'/P_n'$ ($O_n'$ and
$P_n'$ are respectively the ring of integers and prime ideal in $K_n'$).
Then
$$
\xi -g\xi =\pi_n^{-r}\left[ (\sigma_0-g\sigma_0)+(\sigma_1-g\sigma_1)\pi_n
+\cdots +(\sigma_{r-1}-g\sigma_{r-1})\pi_n^{r-1}+\cdots \right]
$$
whence $|\sigma_j-g\sigma_j|_{K_n'}<1$ for $j=0,1,\ldots ,r-1$.

In fact either $\sigma_j\in \widehat{O}_n$, and then $\sigma_j-
g\sigma_j=0$, or $|\sigma_j-g\sigma_j|_{K_n'}=1$ because the
Frobenius automorphism of the unramified extension just permutes
classes from $O_n'/P_n'$ (see Section I-4 in Ref. 16). Therefore all
$\sigma_j$ belong to $\widehat{O}_n$, which means that $\xi \in K_n+O$ as
desired. $\blacksquare$

\section{L\'EVY MEASURE}

{\bf 4.1.} Let us calculate the L\'evy measure of the complement
$S\setminus S_n$.

\begin{teo}
For any $n=1,2,\ldots$
\begin{equation}
\Pi (S\setminus S_n)=(1-q_n^{-1})q_n^{-
ne_n}\frac{q_n^{(ne_n+1)(\alpha /m_n+1)}-q_n^{\alpha
/m_n+1}}{q_n^{\alpha /m_n+1}-1}.
\end{equation}
As $n\to \infty$,
\begin{equation}
\Pi (S\setminus S_n)\sim q_1^{\alpha n}.
\end{equation}
\end{teo}

{\it Proof}. As we saw in Section 2, $\Pi (\K\setminus S)=0$. We
have $S_n\subset S\subset S^{(n)}$, so that
$$
\Pi (S\setminus S_n)=\Pi (S^{(n)}\setminus S_n)-\Pi (S^{(n)}\setminus S)=
\Pi (S^{(n)}\setminus S_n).
$$
The indicator of the set $S^{(n)}\setminus S_n$ is a cylindrical
function (in the sense of Section 2), and we can use (2.6)
obtaining that
\begin{equation}
\Pi (S\setminus S_n)=-\int\limits_{|\eta |_n>1}\|\eta \|^{\alpha
}w_n(\eta )\,d\eta ,
\end{equation}
$$
w_n(\eta )=q_n^{-d_n}\int\limits_{q_n^{d_n-ne_n}<|y|_n\le
q_n^{d_n}}\chi \circ \Tr_{K_n/K_1}(-\eta y)\,dy.
$$

By standard integration formulas (see e.g. Ref. 8, 10)
$$
\int\limits_{|y|_n\le
q_n^{d_n}}\chi \circ \Tr_{K_n/K_1}(-\eta y)\,dy=0,\quad |\eta
|_n>1;
$$
$$
\int\limits_{|y|_n\le
q_n^{d_n-ne_n}}\chi \circ \Tr_{K_n/K_1}(-\eta
y)\,dy=\begin{cases}
q_n^{d_n-ne_n}, & \text{if $|\eta |_n\le q_n^{ne_n}$;}\\
0, & \text{if $|\eta |_n>q_n^{ne_n}$.}
\end{cases}
$$
Substituting into (4.3) we find that
\begin{multline*}
\Pi (S\setminus S_n)=q_n^{-ne_n}\int\limits_{1<|\eta |_n\le
q_n^{ne_n}}|\eta|_n^{\alpha /m_n}d\eta =q_n^{-ne_n}\sum\limits_{j=1}^{ne_n}
\int\limits_{|\eta |_n=q_n^j}|\eta|_n^{\alpha /m_n}d\eta \\
=(1-q_n^{-1})q_n^{-ne_n}\sum\limits_{j=1}^{ne_n}q_n^{j(\alpha
/m_n+1)},
\end{multline*}
which results in (4.2).

Since the sequence (1.1) is strictly increasing, we have $m_n\to
\infty$, thus either $f_n\to \infty$, or $e_n\to \infty$ (or
both). If $f_n\to \infty$, we have $q_n=q_1^{f_n}\to \infty$.
Writing (4.1) as
$$
\Pi (S\setminus S_n)=(1-q_n^{-1})\frac{q_1^{\alpha
/e_n}}{q_1^{\alpha /e_n}-q_n^{-1}}\left( q_1^{\alpha n}-q_1^{-
nm_n}\right)
$$
we come to (4.2) in both possible cases ($e_n\to \infty$ or $e_n$
is constant starting from some value of $n$). If $f_n\equiv f$
for $n\ge n_0$, and $e_n\to \infty$, then we get for $n\ge n_0$
$$
(1-q_{n_0}^{-1})\frac{q_1^{\alpha
/e_n}}{q_1^{\alpha /e_n}-q_{n_0}^{-1}}\longrightarrow 1,\quad
n\to \infty,
$$
and we again come to (4.2). $\quad \blacksquare$

\medskip
{\bf 4.2.} Define a sequence $\{n(j)\}_1^\infty$ setting $n(1)=1$,
$$
n(j+1)=\inf \left\{ n>n(j):\ \Pi (S\setminus S_n)/\Pi (S\setminus
S_{n(j)})\ge 2\right\},\quad j=1,2,\ldots ,
$$
and set $b_n=\left[ \Pi (S\setminus S_{n(j)})\right]^{-1}\log j$
for $n(j)\le n<n(j+1)$.

Evans$^{(3)}$ showed that almost surely
\begin{equation}
\limsup\limits_{n\to \infty}\frac{\pi (n)}{b_n}=1
\end{equation}
where $\pi(n)$ is the first exit time of the process $X_S(t)$ out
of the subgroup $S_n$. Obviously we may substitute for $b_n$ any
sequence $B_n$ with $B_n\sim b_n$, $n\to \infty$. In particular,
if $\alpha >\log_{q_1}2$, we may take (in view of (4.2))
\begin{equation}
B_n=q_1^{-\alpha n}\log n,
\end{equation}
so that
$$
\limsup\limits_{n\to \infty}\frac{\pi (n)}{B_n}=1.
$$

Let $\dim$ and $\Dim$ be the Hausdorff and packing dimensions on
$S$ with respect to the metric defined above (see Ref. 3). Then
$\dim S=\Dim S=1$. It was shown by Evans$^{(3)}$ that for each
$t>0$ almost surely
$$
\dim X([0,t])=\beta',\quad \Dim X([0,t])=\beta'',
$$
where
$$
\beta'=\inf \left\{ \beta :\ \liminf\limits_{n\to
\infty}[M(n)]^{-\beta }Q(n,N)\Pi (S\setminus S_n)=0\right\},
$$
$$
\beta''=\inf \left\{ \beta :\ \limsup\limits_{n\to
\infty}[M(n)]^{-\beta }Q(n,N)\Pi (S\setminus S_n)=0\right\},
$$
where $N$ is an arbitrary natural number. Since by definition
$Q(n,N)\le 1$, the following result is a consequence of Theorem
1.

\begin{cor}
For each $t>0$
$$
\dim X([0,t])=\Dim X([0,t])=0
$$
almost surely.
\end{cor}

\medskip
Unfortunately the results of Ref. 3 do not yield the uniform
dimension results in a similar way. The reason is that
$m_{n+1}\ge 2m_n$ whence
$$
\lim\limits_{n\to \infty }\frac{M(n+1)}{[M(n)]^{1+\eta }}\ne 0
$$
if $\eta >0$ is small enough. This contradicts an a priori
assumption in the uniform dimension study of Ref. 3.

\section{HAUSDORFF MEASURE}

{\bf 5.1.} Let
$$
\tau (n,N)=\meas \left\{ t\in [0,\pi (N)]:\ X(t)\in
S_n\right\},\quad n>N.
$$
It was shown by Evans$^{(3)}$ that
\begin{equation}
[Q(n,N)]^{-1}=\mathbf E[\tau (n,N)]\Pi (S\setminus S_n).
\end{equation}
Moreover (see the proof of Lemma 9 in Ref. 3),
\begin{equation}
\mathbf E[\tau (n,N)]=[M(n)]^{-1}\sum\limits_{\xi \in
\Xi_n}\lambda_n (\xi ),
\end{equation}
where
\begin{equation}
\lambda_n (\xi )=\left[ \Pi (S\setminus S_N)+\int\limits_{S_N}(1-
\chi (\langle x,\xi \rangle ))\Pi (dx)\right]^{-1}.
\end{equation}
Here we identify a class $\xi +O\in \Xi_n$ with an element $\xi$.
Correspondingly, the summation in (5.2) is taken actually over
the set of different classes $\xi +O$.

In order to get an estimate of $\mathbf E[\tau (n,N)]$, we have
to begin with estimating the integral
$$
I_N(\xi )=\int\limits_{S_N}(1-\chi (\langle x,\xi \rangle ))\Pi
(dx),\quad \xi \in \Xi_n,\ n>N.
$$
It is clear that $I_N(\xi )$ is a non-negative real-valued
function.

\begin{lem}
Suppose that the extensions $K_n/K_N$ and $K_N/K_1$ are tamely
ramified. If $\xi \in K_n$, $|\xi |_n=q_n^j$, $Ne_n+1\le j\le
ne_n$, then
\begin{equation}
I_N(\xi )\ge \left( 1-q_N^{-Ne_N}\right) |\xi |_n^{\alpha /m_n}-
\Pi (S\setminus S_N).
\end{equation}
\end{lem}

{\it Proof}. Writing $\chi (\langle x,\xi \rangle ))=\chi (T(\xi
T_n(x))$ we can use (2.6) and come to the expression
$$
I_N(\xi )=-\int\limits_{\zeta \in K_n:\
|\zeta|_n>1}|\zeta|_n^{\alpha /m_n}\Phi (\zeta)\,d\zeta
$$
where
$$
\Phi (\zeta)=q_n^{-d_n}\int\limits_{V_{n,N}}\chi \circ
\Tr_{K_n/K_1}(-\zeta y)\left[ 1-\chi \circ
\Tr_{K_n/K_1}(\xi y)\right]\,dy
$$
and
$$
V_{n,N}=\left\{ y\in K_n:\ |y|_n\le q_n^{d_n},\ \left|
\Tr_{K_n/K_N}(y)\right|_N\le q_N^{d_N-Ne_N}\right\}.
$$
Denoting
$$
\Phi_1(\zeta)=q_n^{-d_n}\int\limits_{V_{n,N}}\chi \circ
\Tr_{K_n/K_1}(-\zeta y)\,dy
$$
we find that
$$
I_N(\xi )=\int\limits_{|\zeta|_n>1}|\zeta|_n^{\alpha /m_n}\left[ \Phi_1
(\zeta -\xi )-\Phi_1(\zeta )\right] \,d\zeta .
$$

Since $V_{n,N}$ is a compact open subgroup in $\Sigma_{n,0}$, we
have (see e.g. (31.7) in Ref. 6)
$$
\Phi_1(\zeta )=\begin{cases}
0, & \text{if $\zeta \notin V_{n,N}^\bot$,}\\
c_{n,N}, & \text{if $\zeta \in V_{n,N}^\bot$,}\end{cases}
$$
where $V_{n,N}^\bot$ is the annihilator of $V_{n,N}$ in the dual
group $\Sigma_{n,0}^*$, and as before,
\begin{equation}
V_{n,N}^\bot =\left\{ \zeta =\zeta' +\zeta'':\ \zeta'\in
K_N,1<|\zeta'|_N\le q_N^{Ne_N};\ \zeta''\in K_n,|\zeta''|_n\le 1\right\},
\end{equation}
$$
c_{n,N}=q_n^{-d_n}\int\limits_{V_{n,N}}dy.
$$

Now
\begin{multline*}
I_N(\xi )=c_{n,N}\left[ \int\limits_{\{|\zeta|_n>1\}\cap \left( \xi
+V_{n,N}^\bot \right)}|\zeta|_n^{\alpha /m_n}d\zeta
-\int\limits_{\{|\zeta|_n>1\}\cap
V_{n,N}^\bot }|\zeta|_n^{\alpha /m_n}d\zeta \right] \\
=c_{n,N}\int\limits_{\{|\zeta|_n>1\}\cap
V_{n,N}^\bot }\left( |\xi +\zeta|_n^{\alpha /m_n}-|\zeta|_n^{\alpha
/m_n}\right) \,d\zeta
\end{multline*}
since by our assumption $|\xi|_n\ge q_n^{Ne_n+1}>1$ while for an
element $\zeta'$ appearing in (5.5) we have
$$
|\zeta'|_n=|\zeta'|_N^{\frac{m_n}{m_N}}\le
q_N^{\frac{Ne_Nm_n}{m_N}}=q_n^{\frac{Ne_Nm_n}{m_Nf_{n,N}}}=q_n^{Ne_n}.
$$

Denote
$$
J_1(\xi )=\int\limits_{\{|\zeta|_n>1\}\cap
V_{n,N}^\bot }|\xi +\zeta|_n^{\alpha /m_n}d\zeta ,\quad
J_2=\int\limits_{\{|\zeta|_n>1\}\cap
V_{n,N}^\bot }|\zeta|_n^{\alpha /m_n}d\zeta .
$$
Writing an element $\zeta'\in K_N,1<|\zeta'|_N\le q_N^{Ne_N}$, as
$$
\zeta'=\pi_N^{-Ne_N}\left( \sigma_0+\sigma_1\pi_N+\cdots
+\sigma_{Ne_N-1}\pi_N^{Ne_N-1}+\cdots \right) ,\quad \sigma_j\in
\widehat{O}_N,
$$
we see that the domain of integration is the union of
non-intersecting ``closed'' unit balls
$B_n^{(Ne_n)}(\sigma_0,\ldots ,\sigma_{Ne_N-1})\subset K_n$
centered at the points $\pi_N^{-Ne_N}\left( \sigma_0+\sigma_1\pi_N+\cdots
+\sigma_{Ne_N-1}\pi_N^{Ne_N-1}\right)$ where at least one of the
elements $\sigma_j$ is different from 0. The total quantity of
such balls equals
$$
\left( q_N-1\right) \left( q_N^{Ne_N-1}+q_N^{Ne_N-2}+\cdots +1\right)
=q_N^{Ne_N}-1,
$$
so that
\begin{equation}
J_1(\xi )=|\xi|_n^{\alpha /m_N}\left( q_N^{Ne_N}-1\right) .
\end{equation}

Similarly, since the equality $|\zeta'|_N=q_N^l$ implies
$|\zeta'|_n=q_n^{le_{n,N}}$, we get
\begin{equation}
J_2=\left( q_N-1\right) q_N^{\alpha /m_N}\frac{q_N^{Ne_N(1+\alpha
/m_N)}-1}{q_N^{1+\alpha /m_N}-1}.
\end{equation}

Next we shall find a lower bound for $c_{n,N}$. It follows from
Lemma 1 that any element $y\in V_{n,N}$ can be written as
$y=y'+y''$ where
$$
|y'|_n\le q_n^{d_n},\quad \Tr_{K_n/K_N}(y')=0,\quad |y''|_n\le
q_n^{d_n-Ne_n}.
$$
We shall find a finite set $F$ of elements $y'$, such that $|y_1'-
y_2'|_n>q_n^{d_n-Ne_n}$ for any different $y_1',y_2'\in F$. Then
$V_{n,N}$ will contain non-intersecting balls centered at $y'\in
F$ with the radii $q_n^{d_n-Ne_n}$, so that
\begin{equation}
c_{n,N}\ge (\card F)q_n^{-Ne_n}.
\end{equation}

Let $x\in K_n$, $q_n^{d_n-Ne_n}<|x|_n\le q_n^{d_n}$. Since the
extension $K_n/K_N$ is tamely ramified, we can choose in $K_N'$
(the inertia subfield of $K_n$) and $K_n$ such prime elements
$\pi_N'$ and $\pi_n$ that $\pi_n^{e_{n,N}}=\pi_N'$. If
$|x|_n=q_n^m$, we can write the canonical representation
$x=\pi_n^{-m}(x_0+x_1\pi_n+\cdots )$, $x_j\in \widehat{O}_N'$ (a
complete set of representatives of residue classes from
$O_N'/P_N'$). Below we assume that $\widehat{O}_N\subset
\widehat{O}_N'$.

It is known (Ref. 16, Chapter 8) that for the element $\Tr_{K_n/K_N'}(x)\in
K_N'$ the inequality
$$
\left| \Tr_{K_n/K_N'}(x)\right|_n\le |x|_n
$$
holds. Thus either $\left| \Tr_{K_n/K_N'}(x)\right|_n<|x|_n$, or
$\left| \Tr_{K_n/K_N'}(x)\right|_n=q_n^m$, and in the latter case
$$
\left| \Tr_{K_n/K_N'}(x)\right|_{K_N'}=q_n^\nu ,\quad \nu \in \mathbb
Z,
$$
so that $m=\nu e_{n,N}$, which implies further that $x=\pi_n^{-
m}x_0+z=\left( \pi'_N\right)^{-\nu}x_0+z$, $|z|_n\le q_n^{m-1}$.
Thus in this case
\begin{equation}
\Tr_{K_n/K_N'}(x)=\left( \pi'_N\right)^{-
\nu}x_0+\Tr_{K_n/K_N'}(z),
\end{equation}
so that
\begin{equation}
\left| x-\Tr_{K_n/K_N'}(x)\right|_n\le q_n^{m-1}.
\end{equation}

Let us consider the element
\begin{equation}
y'=x-g\left[ \Tr_{K_n/K_N'}(x)\right]
\end{equation}
where $g$ is the Frobenius automorphism of the extension
$K_N'/K_N$. We have
$$
\Tr_{K_n/K_N}(y')=\Tr_{K_N'/K_N}\left\{ \Tr_{K_n/K_N'}(x)-g\left[
\Tr_{K_n/K_N'}(x)\right]\right\} =0.
$$
As we saw, $|y'|_n=q_n^m$, except in the case, in which $m=\nu
e_{n,N}$. Now we consider that exceptional case.

Since
$$
y'=\left\{ x-\Tr_{K_n/K_N'}(x)\right\} +\left\{
\Tr_{K_n/K_N'}(x)-g\left[ \Tr_{K_n/K_N'}(x)\right]\right\} ,
$$
and the first summand satisfies (5.10), we have only to study the
second summand which, due to (5.9), is equal, up to a ``small''
term, to $\left( \pi_N'\right)^{-\nu}x_0-g\left[ \left(
\pi_N'\right)^{-\nu}x_0\right]$. By the definition of the
Frobenius automorphism,
$$
\left| \left( \pi_N'\right)^{-\nu}x_0-g\left[ \left(
\pi_N'\right)^{-\nu}x_0\right]\right|_{K_N'}=q_n^\nu ,
$$
or, equivalently,
$$
\left| \left( \pi_N'\right)^{-\nu}x_0-g\left[ \left(
\pi_N'\right)^{-\nu}x_0\right]\right|_n=q_n^m,
$$
unless $\left( \pi_N'\right)^{-\nu}x_0\in K_N$.

Thus we have found that $|y'|_n=q_n^m$ except for the case, in
which $m=\nu e_{n,N}$, $\nu \in \mathbb Z$, and $\left(
\pi_N'\right)^{-\nu}x_0\in K_N$. Since the extension $K_N'/K_N$
is unramified, we may write
$$
\left( \pi_N'\right)^{-\nu}x_0=\pi_N^{-\nu}\left(
\widetilde{x}_0+\widetilde{x}_1\pi_N+\cdots \right) ,\quad
\widetilde{x}_j\in \widehat{O}_N',
$$
so that, due to (5.9), for $x$ we obtain the representation
$$
x=\pi_N^{-\nu}\widetilde{x}_0+\pi_n^{-m+1}\left( \check
x_1+\check x_2\pi_n+\cdots \right) ,\quad \check x_j\in \widehat{O}_N',
$$
with $|y'|_n=q_n^m$ if $\widetilde{x}_0\notin \widehat{O}_N$. A
difference $y'{}^{(1)}-y'{}^{(2)}$ of two elements of the form (5.11)
satisfies the equality $|y'{}^{(1)}-y'{}^{(2)}|_n=q_n^m$ if the
corresponding $\widetilde{x}_0^{(1)},\widetilde{x}_0^{(2)}\in
\widehat{O}_N'$ are such that $\widetilde{x}_0^{(1)}-\widetilde{x}_0^{(2)}
\notin O_N$.

Let us construct a subset $\widehat{O}_N''\subset \widehat{O}_N'$
with the property that if $a,b\in \widehat{O}_N''$, $a\ne b$,
then $a-b\notin O_N$. Denote for brevity $\varkappa
=O_N/P_N$, $\varkappa'=O_N'/P_N'$. Then $\varkappa \subset
\varkappa'$, $\card \varkappa =q_N$, $\card \varkappa' =q_n$.
Consider the quotient group $\varkappa'/\varkappa$ of the
additive groups of the finite fields $\varkappa',\varkappa$. We
have
$$
\card (\varkappa'/\varkappa)=q_n/q_N=q_N^{f_{n,N}-1}.
$$

Let us choose in $\varkappa'$ a complete set of representatives
of residue classes from $\varkappa'/\varkappa$, and then for each of them
(as a class from $O_N'/P_N'$) take a representative from
$\widehat{O}_N'$. Let $\widehat{O}_N''$ be the resulting set. If
$a,b\in \widehat{O}_N''$, and $a-b\in O_N$, then the classes of
$a$ and $b$ in $\varkappa'$ would belong to $\varkappa$ whence
$a=b$ as desired. Note that $\card \widehat{O}_N''=q_N^{f_{n,N}-
1}$.

Now we take as $F$ the set of all elements (5.11) with
\begin{equation}
x=\pi_n^{-d_n}\sum\limits_{\substack{0\le j\le Ne_n-
1\\e_{n,N}\nmid d_n-j}}x_j\pi_n^j+\sum\limits_{\nu :\ d_n-
Ne_n+1\le \nu e_{n,N}\le d_n}\pi_N^{-\nu}\widetilde{x}_\nu
\end{equation}
where $x_j\in \widehat{O}_N$, $\widetilde{x}_\nu \in
\widehat{O}_N''$. Note that $\left| \pi_N^{-
\nu}\right|_n=q_n^{\nu e_{n,N}}$, so that the orders of all
non-zero terms in (5.12) are different.

It follows from our assumptions that the extension $K_n/K_1$ is
tamely ramified (see Ref. 4), that is $d_n=e_n-1$. Hence the
second sum in (5.12) is taken over those $\nu$ for which
$$
\frac{e_n(1-N)}{e_{n,N}}\le \nu \le \frac{e_n}{e_{n,N}}-
\frac{1}{e_{n,N}},
$$
that is
$$
(1-N)e_N\le \nu \le e_N-1.
$$
The quantity of such numbers $\nu$ is $e_N-1-[(1-N)e_N-1]=Ne_N$.
Correspondingly, the quantity of terms in the first sum of (5.12)
is $Ne_n-Ne_N$. Thus
$$
\card F=q_n^{Ne_n-Ne_N}\cdot q_N^{(f_{n,N}-
1)Ne_N}=q_n^{Ne_n}\cdot q_N^{-Ne_N},
$$
and by (5.8)
\begin{equation}
c_{n,N}\ge q_N^{-Ne_N}.
\end{equation}

Comparing (5.6), (5.7), and (5.13) with (4.1) we come to (5.4).
$\quad \blacksquare$

\medskip
Now we can prove (1.2).

\begin{teo}
If all the extensions (1.1) are tamely ramified, then (1.2) holds.
\end{teo}

{\it Proof}. From (5.1), (5.2), (5.3), (4.2), and (3.3) we find
that
\begin{equation}
[Q(n,N)]^{-1}\le C_Nq_1^{\alpha n-nm_n}\sum\limits_{\xi \in
\Xi_n}\left[ \Pi (S\setminus S_N)+I_N(\xi )\right]^{-1},\quad
C_N>0.
\end{equation}
As before, we identify a class $\xi +O$, $\xi \in K_n$,
$|\xi|_n=q_n^j$, $j\ge 1$, with the element
$$
\xi =\pi_n^{-j}\left( \xi_0+\xi_1\pi_n+\cdots +\xi_{j-1}\pi_n^{j-
1}\right) ,\quad \xi_j\in \widehat{O}_n,\xi_0\ne 0.
$$
The number of such elements is $(q_n-1)q_n^{j-1}$.

Let us split the sum in (5.14) into two sums, over $\xi$ with
$|\xi|_n\le q_n^{Ne_n}$, and with $q_n^{Ne_n+1}\le |\xi|_n\le
q_n^{ne_n}$. The first sum is estimated by dropping $I_n(\xi )$;
an upper bound for the second sum is given by Lemma 2. We find
that
\begin{multline*}
[Q(n,N)]^{-1}\le C_Nq_1^{\alpha n-nm_n}\left\{\left[ \Pi
(S\setminus S_N)\right]^{-1}(q_n-
1)\sum\limits_{j=1}^{Ne_n}q_n^{j-1}\right. \\
\left. +(q_n-1)\left( 1-q_N^{-
Ne_N}\right) \sum\limits_{j=Ne_n+1}^{ne_n}q_n^{j-1-j\alpha
/m_n}\right\}\\
\le C_N'q_1^{\alpha n-nm_n}\left[
q_1^{Nm_n}+\frac{q_n^{(ne_n+1)(1-\frac{\alpha }{m_n})}-
q_n^{(Ne_n+1)(1-\frac{\alpha }{m_n})}}{q_n^{1-\frac{\alpha
}{m_n}}-1}\right] \\
\le C_N''q_1^{\alpha n-nm_n}\left( q_1^{Nm_n}+q_1^{nm_n-\alpha
n}\right) \le C_N'''
\end{multline*}
(with positive constants $C_N',C_N'',C_N'''$; we used the fact
that $m_n\to \infty$), which implies (1.2).$\quad \blacksquare$

\medskip
{\bf 5.2.} In accordance with the general definition given in
Ref. 3, a Hausdorff outer measure with respect to a
non-decreasing function $\varphi :\ \{M(n)^{-1}\}_{n=0}^\infty
\cup \{0\}\to [0,\infty )$, such that $\lim\limits_{n\to \infty}
\varphi \left( M(n)^{-1}\right) =\varphi (0)=0$, is defined as
$$
\varphi -m(A)=\liminf\limits_{n\to \infty}\left\{
\sum\limits_i\varphi (\diam R_i)\right\}
$$
where the infimum is taken over all countable collections of
balls $\{R_i\}$ such that a set $A\subset S$ is contained in
$\bigcup\limits_iR_i$, and $\sup\limits_i\diam R_i\le [M(n)]^{-
1}$. Here balls and diameters on $S$ are understood in the sense
of the metric $\Delta$.

In view of Theorem 2, the results by Evans$^{(3)}$ yield the
following construction of the Hausdorff measure for our
situation.

\begin{cor}
In the notation of Section 4 define the function $\varphi$ by
$\varphi (M(n)^{-1})=b_n$ and $\varphi (0)=0$. If all the
extensions (1.1) are tamely ramified, then almost surely
$$
0<\varphi -m(X_S([0,t]))<\infty
$$
for all $t>0$. If $\alpha >\log_{q_1}2$, then the sequence $\{b_n\}$
can be replaced by the sequence $\{B_n\}$ defined by (4.5).
\end{cor}

\section*{REFERENCES}

\begin{enumerate}
\item Albeverio, S., Zhao, X. (2000). On the relation
between different constructions of random walks on $p$-adics.
{\it Markov Processes Relat. Fields} {\bf 6}, 239--255.
\item Evans, S. N. (1989). Local Field Gaussian Measures.
In Cinlar, E. et al. (eds.), {\it Seminar on Stochastic Processes 1988},
Birkh\"{a}user, Boston, pp. 121--160.
\item Evans, S. N. (1989). Local Properties of L\'evy Processes
on a Totally Disconnected Group. {\it J. Theor. Probab.} {\bf 2}, 209--259.
\item Fesenko, I.B. and Vostokov, S.V. (1993). {\it Local
Fields and Their Extensions: A Constructive Approach},
American Mathematical Society, Providence.
\item Haran, S. (1993). Analytic potential theory over the $p$-adics.
{\it Ann. Inst. Fourier} {\bf 43}, 905--944.
\item Hewitt, E. and Ross, K. A. (1970). {\it Abstract
Harmonic Analysis}. Volume 2, Springer, Berlin.
\item Kochubei, A.N. (1992). Parabolic equations over the field of
$p$-adic numbers. {\it Math. USSR Izvestiya} {\bf 39}, 1263--1280.
\item Kochubei, A.N. (1999). Analysis and probability over
infinite extensions of a local field, {\it Potential Anal.} {\bf
10}, 305--325.
\item Kochubei, A.N. (1999). Fractional differentiation
operator over an infinite extension of a local field. In
K\c{a}kol, J. et al. (eds), {\it p-Adic Functional Analysis},
Lect. Notes Pure Appl. Math. Vol. 207, Marcel Dekker, New York,
pp. 167--178.
\item Kochubei, A.N. (2001). {\it Pseudo-Differential Equations
and Stochastics over Non-\linebreak Archimedean Fields}, Marcel Dekker, New York.
\item Lang, S. (1965). {\it Algebra}, Addison-Wesley,
Reading.
\item Narkiewicz, W. (1974). {\it Elementary and Analytic Theory
of Algebraic Numbers}, PWN, Warsaw.
\item Serre, J.-P. (1979). {\it Local Fields}, Springer,
New York.
\item Vladimirov, V.S., Volovich, I.V., and Zelenov, E.I. (1994).
{\it $p$-Adic Analysis and Mathematical Physics}, World Scientific,
Singapore.
\item van der Waerden, B. L. (1971). {\it Algebra I},
Springer, Berlin.
\item Weil, A. (1967). {\it Basic Number Theory}, Springer, Berlin.
\item Yasuda, K. (1996). Additive processes on local fields.
{\it J. Math. Sci. Univ. Tokyo} {\bf 3}, 629--654.
\end{enumerate}
\end{document}